\documentclass[a4paper,11pt]{article}
\usepackage{latexsym}
\usepackage{amssymb}
\usepackage{amsfonts}
\usepackage{amsmath}
\usepackage{indentfirst}
\newtheorem{prop}{Proposition}[section]
\newtheorem{teor}{Theorem}[section]
\newtheorem{cor}{Corollary}[section]
\newcommand{\cvd}{\quad $\blacksquare$}
\newcommand{\cngamma}{\mathcal{C}_{n}^{\Gamma}}

\date{}
\author{Luca Ferrari\thanks{Dipartimento di Matematica ``U.
Dini",
Viale Morgagni 67/A, 50135 Firenze, Italy {\tt
ferrari@math.unifi.it}}\and Renzo Pinzani\thanks{Dipartimento
di Sistemi e Informatica, via Lombroso 6/17, 50135 Firenze, Italy
{\tt pinzani@dsi.unifi.it}}}
\title{Lattices of lattice paths\footnote{This work was partially supported
by MIUR project: \emph{Linguaggi formali e automi: teoria e applicazioni}.}}
\begin{document}
\maketitle

\begin{abstract}
We consider posets of lattice paths (endowed with a natural order) and begin
the study of such structures. We give an algebraic condition to recognize
which ones of these posets are lattices. Next we study the class of Dyck
lattices (i.e., lattices of Dyck paths) and give a recursive construction
for them. The last section is devoted to the presentation of a couple of
open problems.
\end{abstract}

\begin{center}
\end{center}

\section{Introduction}

When a class of objects is introduced in mathematics, one of the
first problems that naturally arises is to count how many objects
there are. Typically, one recognizes an interesting numerical
parameter and tries to enumerate the objects according to it, so
obtaining a sequence of nonnegative integers. This is the main
problem of enumerative combinatorics.

The second step (after mere enumeration) is to look for some ``mathematical"
structure (like, e.g., operations) the class of objects naturally possesses.
As a matter of fact, one of the simplest structure that can be found is an
order relation, as Gian-Carlo Rota suggested in his masterful paper
\cite{Moebius}. By the way, his theory of M\"obius functions of posets
leads to many deep enumerative results, so proving that a better knowledge
of the mathematical structure of a set gives a better insight even on his
combinatorial and enumerative properties.

In this work we deal with classes of lattice paths (like, e.g.,
Dyck paths, Motzkin paths, Schr\"oder paths, Lukasiewicz paths,
etc.) which can be ordered in a completely natural manner. We
postpone the formal definitions of these classes of paths at the
beginning of section \ref{char}. However, in fig. 1 we give some
instances of the paths which will more frequently occur in the
present work. More precisely, we have drawn: \emph{(a)} a Dyck
path, \emph{(b)} a Motzkin path, \emph{(c)} a Schr\"oder path and
\emph{(d)} a Lukasiewicz path.

\begin{center}
\setlength{\unitlength}{1.8mm}
\begin{picture}(60,30)
\put(0,19){\circle*{0.3}} \put(2,21){\circle*{0.3}}
\put(4,23){\circle*{0.3}} \put(6,25){\circle*{0.3}}
\put(8,23){\circle*{0.3}} \put(10,25){\circle*{0.3}}
\put(12,23){\circle*{0.3}} \put(14,21){\circle*{0.3}}
\put(16,23){\circle*{0.3}} \put(18,21){\circle*{0.3}}
\put(20,19){\circle*{0.3}} \put(22,21){\circle*{0.3}}
\put(24,19){\circle*{0.3}} \put(0,19) {\line(1,1){2}}
\put(2,21){\line(1,1){2}} \put(4,23){\line(1,1){2}}
\put(6,25){\line(1,-1){2}} \put(8,23){\line(1,1){2}}
\put(10,25){\line(1,-1){2}} \put(12,23){\line(1,-1){2}}
\put(14,21){\line(1,1){2}} \put(16,23){\line(1,-1){2}}
\put(18,21){\line(1,-1){2}} \put(20,19){\line(1,1){2}}
\put(22,21){\line(1,-1){2}} \put(11,17){\makebox(0,0){\emph{(a)}}}

\put(30,19){\circle*{0.3}} \put(32,21){\circle*{0.3}}
\put(34,23){\circle*{0.3}} \put(36,23){\circle*{0.3}}
\put(38,25){\circle*{0.3}} \put(40,23){\circle*{0.3}}
\put(42,21){\circle*{0.3}} \put(44,21){\circle*{0.3}}
\put(46,21){\circle*{0.3}} \put(48,19){\circle*{0.3}}
\put(50,21){\circle*{0.3}} \put(52,19){\circle*{0.3}}
\put(54,19){\circle*{0.3}} \put(30,19){\line(1,1){2}}
\put(32,21){\line(1,1){2}} \put(34,23){\line(1,0){2}}
\put(36,23){\line(1,1){2}} \put(38,25){\line(1,-1){2}}
\put(40,23){\line(1,-1){2}} \put(42,21){\line(1,0){2}}
\put(44,21){\line(1,0){2}} \put(46,21){\line(1,-1){2}}
\put(48,19){\line(1,1){2}} \put(50,21){\line(1,-1){2}}
\put(52,19){\line(1,0){2}} \put(41,17){\makebox(0,0){\emph{(b)}}}

\put(0,4){\circle*{0.3}} \put(2,6){\circle*{0.3}}
\put(6,6){\circle*{0.3}} \put(8,8){\circle*{0.3}}
\put(10,10){\circle*{0.3}} \put(14,10){\circle*{0.3}}
\put(16,8){\circle*{0.3}} \put(20,8){\circle*{0.3}}
\put(22,6){\circle*{0.3}} \put(24,4){\circle*{0.3}}
\put(0,4){\line(1,1){2}} \put(2,6){\line(2,0){2}}
\put(4,6){\line(2,0){2}} \put(6,6){\line(1,1){2}}
\put(8,8){\line(1,1){2}} \put(10,10){\line(2,0){2}}
\put(12,10){\line(2,0){2}} \put(14,10){\line(1,-1){2}}
\put(16,8){\line(2,0){2}} \put(18,8){\line(2,0){2}}
\put(20,8){\line(1,-1){2}} \put(22,6){\line(1,-1){2}}
\put(11,2){\makebox(0,0){\emph{(c)}}}

\put(30,4){\circle*{0.3}} \put(32,6){\circle*{0.3}}
\put(34,10){\circle*{0.3}} \put(36,8){\circle*{0.3}}
\put(38,14){\circle*{0.3}} \put(40,12){\circle*{0.3}}
\put(42,10){\circle*{0.3}} \put(44,8){\circle*{0.3}}
\put(46,6){\circle*{0.3}} \put(48,4){\circle*{0.3}}
\put(50,8){\circle*{0.3}} \put(52,6){\circle*{0.3}}
\put(54,4){\circle*{0.3}} \put(30,4){\line(1,1){2}}
\put(32,6){\line(1,2){2}} \put(34,10){\line(1,-1){2}}
\put(36,8){\line(1,3){2}} \put(38,14){\line(1,-1){2}}
\put(40,12){\line(1,-1){2}} \put(42,10){\line(1,-1){2}}
\put(44,8){\line(1,-1){2}} \put(46,6){\line(1,-1){2}}
\put(48,4){\line(1,2){2}} \put(50,8){\line(1,-1){2}}
\put(52,6){\line(1,-1){2}} \put(41,2){\makebox(0,0){\emph{(d)}}}

\put(30,-2){\makebox(0,0){Fig.1}}
\end{picture}
\end{center}

\bigskip

\bigskip

Using a standard vocabulary, we say that Dyck paths only use rise
(or $(1,1)$) steps and fall (or $(1,-1)$) steps, whereas Motzkin
and Schr\"oder paths also use horizontal steps of length 1 (or
$(1,0)$ steps) and of length 2 (or $(2,0)$ steps), respectively.
The case of Lukasiewicz paths is quite different, since they use
an infinite set of possible steps, namely rise step of any type
$(1,k)$ and (simple) fall steps of type $(1,-1)$.

\bigskip

Our goal is to begin a systematic study of the posets of paths
arising in this way. A similar point of view has been undertaken
in \cite{delannoy}, where the authors study the posets arising
from Delannoy paths; however, the present work deals with
essentially different classes of paths. The only general problem
we tackle here is to determine in which cases a poset of paths is
a lattice, and we propose a possible solution to this problem. We
also point out that the first works in this direction have been
done by Narayana (see, for example, \cite{nara}): in \cite{nara2}
it is proved a very interesting result, described at the end of
section \ref{char}. Next, we focus on the study of a single type
of paths, namely Dyck paths, and we provide an explicit
(recursive) construction for the lattices of Dyck paths. It is
then easy to see that a suitable, slight modification of such a
construction can be successfully applied also to lattices of
Schr\"oder paths. We remark that in \cite{partizioni, partbaseb}
some similar problems are studied: the authors provide the
construction of the lattice of partitions of a given integer $n$
(with the dominance order) starting from the knowledge of the
lattice of partitions of $n-1$. However, the methods used in the
present work are completely different from those employed in
\cite{partizioni, partbaseb}.

Quite surprisingly, our study sheds new light on a method of
enumeration, usually called \emph{ECO method} (ECO stands for
``Enumeration of Combinatorial Objects"), which has proved
fruitful in many problems of enumeration. A rough description of
this method is given in section \ref{dyck} (mainly with the help
of the example of Dyck paths). What is important to point out here
is that such a method provides (among other things) a partition of
the objects of size $n$ of a given set of combinatorial objects
(when a suitable definition of size is provided); this fact is
obviously crucial in trying to enumerate the objects of the class
according to the size. What is typical of the ECO method when
applied to the construction of a class of paths (Dyck, Motzkin,
Schr\"oder, and so on) is that the equivalence classes of the
above mentioned partition have a nice order structure: they are
\emph{chains} (i.e., totally ordered sets of paths) with respect
to the natural order of paths introduced in this paper. So the ECO
method, if applicable, provides a chain partition of the lattice
of paths under consideration.

\bigskip

Throughout the whole paper $\mathbf{N}$ will denote the set of
nonnegative integers and $\mathbf{Z}$ the set of all integers. For
$x,y\in \mathbf{Z}$, the expression $[x,y]$ will always indicate
the interval of \emph{integers} between $x$ and $y$.

Given a lattice $\mathbf{L}$, a \emph{filter} of $\mathbf{L}$ is a
subset $F$ such that, if $x,y\in F$, then $x\land y\in F$ (that
is, $F$ is closed for meet) and such that, if $x\in F$ and $x\leq
y$, then $y\in F$. Dually, an \emph{ideal} of $\mathbf{L}$ is a
subset $I$ such that, if $x,y\in I$, then $x\lor y\in I$ (that is,
$F$ is closed for join) and such that, if $x\in I$ and $y\leq x$,
then $y\in I$. If $\mathbf{L}$ has minimum $\mathbf{0}$, $L$ is
said to be \emph{ranked} when there exists a map
$r:\mathbf{L}\longrightarrow \mathbf{N}$ such that
$r(\mathbf{0})=0$ and, if $y$ covers $x$ (i.e., $x<y$ and there is
no element $z$ such that $x<y<z$), then $r(y)=r(x)+1$. The
function $r$ is called the \emph{rank function} of $L$. For any
other poset and lattice concept, we refer to the texts
\cite{CD,priestley}. In \cite{encomb} a whole chapter is devoted
to the study of posets, especially in connection with enumerative
and algebraic combinatorics.

Our last remark concerns the use of the term ``isomorphism" (and similar ones):
in this paper it has to be considered a synonym of \emph{lattice isomorphism},
i.e. a bijective map between two lattices which is both join-preserving and
meet-preserving.

\section{Definitions and a first characterization}\label{char}

Let $\Gamma$ be a finite subset of $\mathbf{Z}$. We call
\emph{$\Gamma$-path of length $n$} every element of the set:
\begin{eqnarray}
\cngamma&=&\{ f:[0,n]\longrightarrow
\mathbf{N}\; |\; f(0)=f(n)=0; \nonumber
\\ & &f(k+1)-f(k)\in \Gamma, \forall k<n\}.
\end{eqnarray}

The elements of the set $\mathcal{C}^{\Gamma}=\bigcup
_{n}\cngamma$ are called \emph{$\Gamma$-paths}. The elements of
the set $\Gamma$ are called \emph{steps}. This definition is
completely equivalent to the usual definition of a lattice path
starting from the origin, ending on the $x$-axis and using steps
of a prescribed type. Observe that this definition allows us to
consider only lattice paths whose steps have length 1 (that is, of
the form $(1,k)$); so, for example, Schr\"oder paths do not fall
within our definition (since horizontal steps are of type
$(2,0)$). The following examples can be immediately checked by the
reader:

\begin{itemize}
\item[$\bullet$)] $\{ -1,1\}$-paths are the ordinary Dyck paths;
\item[$\bullet$)] $\{ -1,0,1\}$-paths are the ordinary Motzkin
paths.
\end{itemize}

It is possible to introduce a natural order between the elements
of each set $\cngamma$, by defining $f\leq g$ in $\cngamma$
whenever $f(i)\leq g(i)$, for any $i\leq n$. The fact that
$[\cngamma ; \leq ]$ is a poset is immediate. Instead, it could be
interesting to wonder whether such an order is a lattice order, that
is: for which $\Gamma \subseteq \mathbf{Z}$ the poset
$[\cngamma ; \leq ]$ is a lattice? In this section we will be concerned
precisely with this problem.

\bigskip

A first remark to be done is that, among the various lattices of
$\Gamma$-paths, there are some cases which are particularly nice. Indeed, for
some choices of $\Gamma$, the meet and join
operations induced on $\cngamma$ by the order introduced above are
defined pointwise, i.e.:
\begin{eqnarray}\label{coord}
(f\lor g)(i)=f(i)\lor g(i), \nonumber
\\ (f\land g)(i)=f(i)\land g(i),
\end{eqnarray}

\bigskip

\bigskip

\begin{center}
\setlength{\unitlength}{0.8mm}
\begin{picture}(120,100)
\put(50,20){\circle*{1}} \put(55,20){\circle*{1}}
\put(60,20){\circle*{1}} \put(65,20){\circle*{1}}
\put(70,20){\circle*{1}} \put(50,20){\line(1,0){5}}
\put(55,20){\line(1,0){5}} \put(60,20){\line(1,0){5}}
\put(65,20){\line(1,0){5}} \put(50,15){\makebox(0,0){{\bf 0}}}

\put(0,40){\circle*{1}} \put(5,50){\circle*{1}}
\put(10,45){\circle*{1}} \put(15,40){\circle*{1}}
\put(20,40){\circle*{1}} \put(0,40){\line(1,2){5}}
\put(5,50){\line(1,-1){5}} \put(10,45){\line(1,-1){5}}
\put(15,40){\line(1,0){5}} \put(0,35){\makebox(0,0){$a$}}

\put(0,70){\circle*{1}} \put(5,80){\circle*{1}}
\put(10,75){\circle*{1}} \put(15,75){\circle*{1}}
\put(20,70){\circle*{1}} \put(0,70){\line(1,2){5}}
\put(5,80){\line(1,-1){5}} \put(10,75){\line(1,0){5}}
\put(15,75){\line(1,-1){5}} \put(0,65){\makebox(0,0){$c$}}

\put(100,40){\circle*{1}} \put(105,40){\circle*{1}}
\put(110,50){\circle*{1}} \put(115,45){\circle*{1}}
\put(120,40){\circle*{1}} \put(100,40){\line(1,0){5}}
\put(105,40){\line(1,2){5}} \put(110,50){\line(1,-1){5}}
\put(115,45){\line(1,-1){5}} \put(100,35){\makebox(0,0){$b$}}

\put(50,90){\circle*{1}} \put(55,100){\circle*{1}}
\put(60,100){\circle*{1}} \put(65,95){\circle*{1}}
\put(70,90){\circle*{1}} \put(50,90){\line(1,2){5}}
\put(55,100){\line(1,0){5}} \put(60,100){\line(1,-1){5}}
\put(65,95){\line(1,-1){5}} \put(50,85){\makebox(0,0){{\bf 1}}}

\put(45,25){\line(-2,1){20}}

\put(75,25){\line(2,1){20}}

\put(10,55){\line(0,1){10}}

\put(25,75){\line(2,1){20}}

\put(96,55){\line(-2,3){20}}

\put(60,5){\makebox(0,0){Fig.2}}
\end{picture}
\end{center}

\noindent where the join and meet in the r.h.s. are computed in
$\mathbf{N}$. Therefore, the lattices arising in such cases are
distributive. Unfortunately, this is not true in general, since
for some choices of $\Gamma$ it happens that $\cngamma$ is
actually a lattice but the operations are not defined as in
(\ref{coord}). The following example will clarify this statement.

\bigskip

{\bf \emph{Example.}} Consider $\Gamma =\{ -1,0,2\}$. The poset
$\mathcal{C}_{4}^{\Gamma}$ contains 5 paths and its Hasse diagram
is depicted in fig. 2. As it is well-known, it is a lattice which
is not distributive (it is not even modular), so lattice
operations on such paths cannot be defined pointwise as in
(\ref{coord}). Indeed, consider the two paths $b$ and $c$: their
meet in $\mathcal{C}_{4}^{\Gamma}$ is the minimum $\mathbf{0}$ of
the lattice, which is not the coordinatewise meet of the two
paths.




\bigskip

The next example shows that there are also cases in which $\cngamma$ is not
a lattice.

\bigskip

{\bf \emph{Example.}} Take $\Gamma =\{ -1,1,2\}$, and consider the
poset $\mathcal{C}_{5}^{\Gamma}$. Its Hasse diagram is the
following:

\begin{center}
\setlength{\unitlength}{0.5mm}
\begin{picture}(60,60)
\put(20,20){\circle*{2}} \put(40,20){\circle*{2}}
\put(20,40){\circle*{2}} \put(40,20){\circle*{2}}
\put(40,40){\circle*{2}} \put(30,50){\circle*{2}}
\put(20,20){\line(0,1){20}} \put(20,20){\line(1,1){20}}
\put(40,20){\line(0,1){20}} \put(20,40){\line(1,1){10}}
\put(40,40){\line(-1,1){10}} \put(30,0){\makebox(0,0){Fig.3}}
\end{picture}
\end{center}
\bigskip
which is not the Hasse diagram of a lattice.

\bigskip

Our next theorem will give an answer to the above problem in the case of
lattices with coordinatewise meet and join. To state and prove our result
we need to introduce some notation. We set $\gamma_{+}=\max \Gamma$,
$\gamma_{-}=\min \Gamma$; it is clear that $\gamma_{+}\geq 0$ and
$\gamma_{-}\leq 0$ (that is, $\Gamma$ must contain both up and down steps,
otherwise the set of $\Gamma$-paths would be empty). We call
\emph{diameter} of $\Gamma$ the difference
between the maximum and the minimum of $\Gamma$, that is
$\overline{\gamma}=\textnormal{diam}\; \Gamma =\gamma_{+}-\gamma_{-}$.
Finally, we define the set $\Delta_{n}^{\Gamma}\subseteq \mathbf{Z}$
as follows (for any $n\in \mathbf{N}$):
\begin{equation}
\Delta_{n}^{\Gamma}=\left\{ \sum_{i=1}^{n}(x_{i}-y_{i})\;
\arrowvert \; x_{i},y_{i}\in \Gamma ,x_{i}\neq y_{j}, \forall i,j
\right\},
\end{equation}
and we set $\Delta ^{\Gamma}=\bigcup_{n\in \mathbf{N}}\Delta _{h}^{\Gamma}$.

\bigskip

{\bf \emph{Example.}} Take $\Gamma =\{ -1,2\}$. Then $\Delta _n
^{\Gamma}$ contains all the sums having $n$ summands, each of the
form $2-(-1)=3$ or $-1-2=-3$. If we let $n$ running over
$\mathbf{N}$, then we obtain $\Delta ^{\Gamma}=3\mathbf{Z}$ (set
of all integers divisible by 3).

\bigskip

Now we are ready to state our main result concerning the algebraic
characterization of a particular class of lattices of lattice paths.

\begin{teor}\label{coor}
$\cngamma$ is a (finite) distributive lattice with respect to
coordinatewise meet and join if and only if
\begin{equation}\label{lattice}
(\Delta ^{\Gamma}+\Gamma)\cap [\gamma_{-},\gamma_{+}]\subseteq \Gamma.
\end{equation}
\end{teor}

\emph{Proof.} It is clear that all the axioms of distributive
lattices are satisfied by $\cngamma$. The only difficult thing to
check in order to prove that $\cngamma$ is a lattice is that it is
closed with respect to the coordinatewise join and meet operations
described above, i.e., if $f,g\in \cngamma$, then $f\land g,f\lor
g\in \cngamma$. So take $f,g\in \cngamma$, and suppose that, for a
given $k\in \mathbf{N}$, $f(k)<g(k)$ and $f(k+1)>g(k+1)$ (this is
the only difficult case to study).

For suitable $\gamma_{1},\gamma_{2}\in \Gamma$, we have
\begin{eqnarray*}
f(k+1)=f(k)+\gamma_{1},
\\ g(k+1)=g(k)+\gamma_{2}.
\end{eqnarray*}

In order that the path $f\lor g$ belongs to $\Gamma$, we must
impose that the step $f(k+1)-g(k)$ is contained in $\Gamma$. Now
$f(k+1)-g(k)=(f(k)-g(k))+\gamma_{1}$. Both $f(k)$ and $g(k)$ can
be expressed as a sum of $k$ elements of $\Gamma$; more precisely,
we can set:
\begin{displaymath}
f(k)=\sum_{i=1}^{k}\tilde{x}_{i},\qquad g(k)=\sum_{i=1}^{k}\tilde{y}_{i},
\end{displaymath}
where $\tilde{x}_{i},\tilde{y}_{i}\in \Gamma$. Thus
we can express the difference $f(k)-g(k)$ as a sum of differences between
elements of $\Gamma$; therefore, possibly deleting pairs of equal steps
appearing both in the expressions of $f(k)$ and $g(k)$, for a suitable
$h\in \mathbf{N}$, we have:
\begin{displaymath}
f(k)-g(k)=\sum_{i=1}^{h}(x_{i}-y_{i}),
\end{displaymath}
where $x_{i}\neq y_{j}$, for every $i,j\leq h$ ($x_{i},y_{j}\in \Gamma$).

Now we are very close to completing our proof. Imposing that
$f(k+1)-g(k)=(f(k)-g(k))+\gamma_{1}\in \Gamma$ leads in general to the
condition
\begin{equation}\label{sup}
((\Delta ^{\Gamma}\cap \; \; ]-\overline{\gamma},0])+\Gamma )\cap
[\gamma_{-},\gamma_{+}]\subseteq \Gamma,
\end{equation}
(recall the definition of the set $\Delta ^{\Gamma}$ given above). The
reasons for which we have considered the intersection $\Delta ^{\Gamma}
\cap \; \; ]-\overline{\gamma},0]$ are the following:
\begin{itemize}
\item[$\bullet$)] $f(k)-g(k)$ must be negative, since we have supposed
that $f(k)<g(k)$;
\item[$\bullet$)] the paths $f$ and $g$ must cross, so the difference
between $f(k)$ and $g(k)$ has to allow the crossing: this is why such
difference must be greater than $-\overline{\gamma}$.
\end{itemize}

Also the fact that we consider the intersection with
$[\gamma_{-},\gamma_{+}]$ needs an explanation. Since we are supposing
there is a crossing, necessarily we have $f(k+1)\geq g(k+1)$ and
$f(k)\leq g(k)$, whence:
\begin{displaymath}
\gamma_{-}\leq g(k+1)-g(k)\leq f(k+1)-g(k)\leq f(k+1)-f(k)\leq
\gamma_{+}.
\end{displaymath}

Thus, condition (\ref{sup}) is necessary for the existence of the
join $f\lor g$ in $\cngamma$. Analogously, it can be shown that a
similar condition must be verified for the existence of the meet
$f\land g$, and precisely:
\begin{equation}
((\Delta ^{\Gamma}\cap [0,\overline{\gamma}[\; )+\Gamma )\cap
[\gamma_{-},\gamma_{+}]\subseteq \Gamma.
\end{equation}

Putting things together, we have proved that, if $\cngamma$ is a
lattice, then:
\begin{equation}
(\Delta ^{\Gamma}+\Gamma)\cap [\gamma_{-},\gamma_{+}]\subseteq \Gamma.
\end{equation}

Conversely, suppose that condition (\ref{lattice}) is satisfied. Then,
if the $\Gamma$-paths $f$ and $g$ cross between $k$ and $k+1$, the fact
that $f(k+1)-g(k)$ and $g(k+1)-f(k)$ both belong to $\Gamma$ is ensured
by the hypothesis, and we are done. \cvd

\begin{cor}\label{1} If $\Gamma$ is an interval (i.e.
$\Gamma=[\gamma_{1},\gamma_{2}]$), then $\cngamma$ is a lattice.
\end{cor}

\emph{Proof.} Indeed, in this case
$\Gamma=[\gamma_{-},\gamma_{+}]$, and so we get
$(\Delta^{\Gamma}+\Gamma )\cap \Gamma \subseteq \Gamma$. \cvd

\begin{cor}\label{2} If $\Gamma =\{ -b,a\}$ ($a,b\in \mathbf{N}$),
then $\cngamma$ is a lattice.
\end{cor}

\emph{Proof.} First observe that, in the assumed hypothesis, we have:
\begin{displaymath}
\Delta^{\Gamma}=\{ (a+b)z\; |\; z\in \mathbf{Z}\} ,
\end{displaymath}
and so
\begin{displaymath}
\Delta^{\Gamma}+\Gamma =\{ (a+b)z+a,(a+b)z-b\; |\; z\in \mathbf{Z}\} .
\end{displaymath}

Now consider the set $(\Delta^{\Gamma}+\Gamma )\cap [-b,a]$, i.e. the
elements of $\Delta^{\Gamma}+\Gamma$ contained in the integer interval
$[-b,a]$. If $c$ is such an element then two cases are possible.
\begin{itemize}
\item[i)] $c=(a+b)z+a$: then $z$ must necessarily be of the form
\begin{displaymath}
z=\frac{c-a}{a+b}\in \mathbf{Z};
\end{displaymath}
the r. h. s. is an increasing function of $c$, moreover $c=-b$
implies that $z=-1$ and $c=a$ implies that $z=0$. Since $z\in
\mathbf{Z}$, this means that no value of $c$ in $[-b,a]$ is allowed
other than $c=-b$ and $c=a$, and so $c\in \{-b,a\} =\Gamma$.
\item[ii)] $c=(a+b)z-b$: a similar argument leads to the same
conclusion, i.e. $c\in \{-b,a\} =\Gamma$. \cvd
\end{itemize}

\bigskip

\underline{\emph{Remark}}. Each of the previous results holds also
if $\Gamma$ is infinite (on the right, on the left, on both
sides). In this case, we have $\gamma_+ =+\infty$ and/or $\gamma_-
=-\infty$; to obtain the corresponding statements, one has to
replace the interval $[\gamma_{-},\gamma_{+}]$ in the text of
theorem \ref{coor} with $[\gamma_- ,+\infty[$, $]-\infty ,\gamma_+
]$ or $]-\infty ,+\infty[$, depending on the cases.

\bigskip

{\bf \emph{Examples.}}
\begin{enumerate}
\item Motzkin paths of a given length $n$ have a lattice
structure. Indeed, as we have seen before, Motzkin paths
correspond to $\Gamma$-paths for $\Gamma =\{ -1,0,1\}$, so in this
case $\Gamma$ is an integer interval and we can apply corollary
\ref{1}. \item Dyck paths are precisely $\{ -1,1\}$-paths, so,
thanks to corollary \ref{2}, they constitute a lattice for every
even length $2n$. Analogously, also $\{2,-1\}$-paths of any given
length (necessarily of the form $3n$) constitute a lattice. \item
Recalling the preceding remark, we can consider also paths with an
infinite set of steps as, for example, the so-called Lukasiewicz
paths (see, for instance, \cite{BF}). According to our vocabulary,
they are $\Gamma$-paths for $\Gamma =\{ -1\} \cup \mathbf{N}$. By
translating corollary \ref{1} into its corresponding version for
the infinite set $\Gamma$, we can conclude that Lukasiewicz paths
of any given length form a lattice.
\end{enumerate}

\underline{\emph{Remark}}. As we have said at the beginning, in
this paper we will be concerned with paths using exclusively
single steps (i.e., of length 1). Anyway, any class of paths of
the same length can be endowed with the natural order we have
defined, therefore it makes sense to ask whether such an order is
a lattice order. For example, it is natural to wonder whether
Schr\"oder paths of a given length possess a lattice structure. In
this case, it is possible to show that the answer is positive,
even if we prefer not to go into details. The reason for which
everything works is that, if $f,g$ are Schr\"oder paths (of the
same length) and $g(k)-f(k)=1$, it means that either $f$ or $g$
possesses a double horizontal step having its middle point in $k$,
so $f$ and $g$ can not cross at $k+1$.

\bigskip

In closing this section, we recall a work of Narayana and Fulton \cite{nara2},
in which it is shown that the set of the (minimal) lattice paths
from the origin to $(m,n)$ endowed with the \emph{dominance order} is a
distributive lattice. Setting $m=n$ and translating into our language, this
means that the set of Grand-Dyck paths (see \cite{2bii} for the definition)
with the order we have defined is a distributive lattice. To the best of our
knowledge, this is the first result towards an order-theoretic
investigation of classes of lattice paths.

\section{Dyck lattices}\label{dyck}

In the previous section we have shown that a lot of classes of paths often
used in combinatorics can be endowed with a lattice structure in a
natural way. In this section we try to examine in details the lattice
structure of a very well-known class of paths, namely Dyck paths. For a very
detailed survey on Dyck paths and their enumeration, see \cite{D}.

We have just proved that, for any given $n\in \mathbf{N}$, Dyck
paths of length $2n$ constitute a finite distributive lattice: we
will denote it by $\mathcal{D}_{n}$. Throughout this whole section
we will always denote the minimum and the maximum of
$\mathcal{D}_{n}$ by $\mathbf{0}$ and $\mathbf{1}$, respectively
(independently from $n$). In fig. 4 we have drawn the Hasse
diagrams of $\mathcal{D}_{n}$ for small values of $n$.

\bigskip

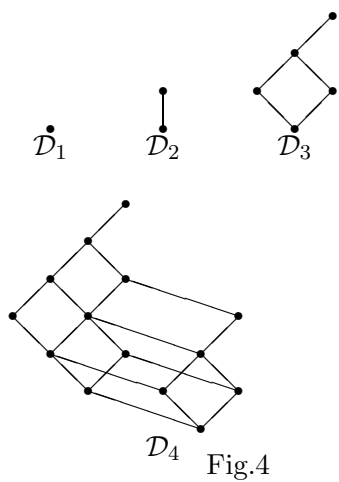
\begin{center}
\setlength{\unitlength}{0.5mm}
\begin{picture}(120,120)
\put(10,90){\circle*{2}} \put(40,90){\circle*{2}}
\put(40,100){\circle*{2}} \put(75,90){\circle*{2}}
\put(65,100){\circle*{2}} \put(85,100){\circle*{2}}
\put(75,110){\circle*{2}} \put(85,120){\circle*{2}}
\put(50,10){\circle*{2}} \put(40,20){\circle*{2}}
\put(60,20){\circle*{2}} \put(50,30){\circle*{2}}
\put(60,40){\circle*{2}} \put(20,20){\circle*{2}}
\put(10,30){\circle*{2}} \put(30,30){\circle*{2}}
\put(20,40){\circle*{2}} \put(30,50){\circle*{2}} \put(
0,40){\circle*{2}} \put(10,50){\circle*{2}}
\put(20,60){\circle*{2}} \put(30,70){\circle*{2}}
\put(40,90){\line(0,1){10}} \put(75,90){\line(1,1){10}}
\put(65,100){\line(1,1){10}} \put(75,110){\line(1,1){10}}
\put(75,90){\line(-1,1){10}} \put(85,100){\line(-1,1){10}}
\put(50,10){\line(1,1){10}} \put(40,20){\line(1,1){10}}
\put(50,30){\line(1,1){10}} \put(20,20){\line(1,1){10}}
\put(10,30){\line(1,1){10}} \put(20,40){\line(1,1){10}} \put(
0,40){\line(1,1){10}} \put(10,50){\line(1,1){10}}
\put(20,60){\line(1,1){10}} \put(50,10){\line(-1,1){10}}
\put(60,20){\line(-1,1){10}} \put(20,20){\line(-1,1){10}}
\put(30,30){\line(-1,1){10}} \put(10,30){\line(-1,1){10}}
\put(20,40){\line(-1,1){10}} \put(30,50){\line(-1,1){10}}
\put(50,10){\line(-3,1){30}} \put(40,20){\line(-3,1){30}}
\put(60,20){\line(-3,1){30}} \put(50,30){\line(-3,1){30}}
\put(60,40){\line(-3,1){30}}
\put(10,85){\makebox(0,0){$\mathcal{D}_{1}$}}
\put(40,85){\makebox(0,0){$\mathcal{D}_{2}$}}
\put(75,85){\makebox(0,0){$\mathcal{D}_{3}$}}
\put(40,5){\makebox(0,0){$\mathcal{D}_{4}$}}
\put(60,0){\makebox(0,0){Fig.4}}
\end{picture}
\end{center}

\bigskip

Our goal is to give a description of the shape of the lattices
$\mathcal{D}_{n}$ and then to provide an efficient way for constructing
$\mathcal{D}_{n}$ starting from the knowledge of $\mathcal{D}_{n-1}$,
in the same spirit of \cite{partizioni}, where the same problem is solved
(in a completely different way) for lattices of integer partitions.

\bigskip

The first thing we observe is that $\mathcal{D}_{n}$ has a
\emph{rank function}: this follows immediately from the fact that
it is a distributive lattice (see \cite{priestley}). The
determination of the rank of an element of $\mathcal{D}_{n}$ is
the object of the next proposition.

\begin{prop}\label{rango} If we denote by $r_{n}$ the rank function of
$\mathcal{D}_{n}$, then, given $f\in \mathcal{D}_{n}$, we have:
\begin{equation}
r_{n}(f)=\frac{\alpha (f)-n}{2},
\end{equation}
where $\alpha (f)$ denotes the \emph{area} of $f$, that is, by definition,
the sum of the values of $f$:
\begin{displaymath}
\alpha (f)=\sum_{k=1}^{n}f(k).
\end{displaymath}
\end{prop}

\emph{Proof.} The minimum $\mathbf{0}\in \mathcal{D}_{n}$ is the
function defined as follows (for $k\leq 2n$):
\begin{displaymath}
\mathbf{0}(k)=
\left\{ \begin{array}{ll}
0, & \textnormal{if $k$ is even;}
\\ 1, & \textnormal{if $k$ is odd.}\end{array} \right.
\end{displaymath}

Thus we have $r_{n}(\mathbf{0})=0=\frac{n-n}{2}=\frac{\alpha
(\mathbf{0})-n}{2}$. Now, by induction, suppose that
$r_{n}(f)=\rho=\frac{\alpha (f)-n}{2}$ for a given $f\in
\mathcal{D}_{n}$, and take $g\in \mathcal{D}_{n}$ such that
$f\prec g$ (this is a standard notation to mean that $f$ is
covered by $g$). Then, there exists precisely one $k$ such that
$2\leq k\leq 2n-2$, $f(h)=g(h)$ for every $h\neq k$ and
$g(k)=f(k)+2$. Therefore, by the definition of area given above,
$\alpha (g)=\alpha (f)+2$ and, by the definition of rank function,
\begin{displaymath}
r_{n}(g)=\rho +1=\frac{\alpha (f)-n}{2}+1=\frac{\alpha
(f)+2-n}{2}=\frac{\alpha (g)-n}{2}. \textnormal{\cvd}
\end{displaymath}

\bigskip

\underline{\emph{Remark}}. Thanks to the above proposition, we can
assert that two Dyck paths of the same length have the same rank
if and only if they have the same area.

\bigskip

There is a very classical and natural way of investigating the
structure of finite distributive lattices. A famous result by
Birkhoff asserts that each finite distributive lattice is
isomorphic to the lattice of the down-sets (or order ideals) of
the poset of its join-irreducible elements. In the case of Dyck
lattices, join-irreducible elements are quite easy to describe:
they are precisely those paths having exactly one hill of height
$>1$ (where by ``hill" of height $h$ we mean a sequence of $h$
rise steps followed by a sequence of $h$ fall steps). A careful
analysis of the situation shows that the poset of the
join-irreducible elements of $\mathcal{D}_{n}$ is precisely the
poset of the intervals of the $n$-element chain. For $n=4$, for
instance, the Hasse diagram of this poset is the following:

\begin{center}
\setlength{\unitlength}{0.5mm}
\begin{picture}(60,50)
\put(10,20){\circle*{2}} \put(30,20){\circle*{2}}
\put(50,20){\circle*{2}} \put(20,30){\circle*{2}}
\put(40,30){\circle*{2}} \put(30,40){\circle*{2}}
\put(10,20){\line(1,1){10}} \put(30,20){\line(1,1){10}}
\put(30,20){\line(-1,1){10}} \put(50,20){\line(-1,1){10}}
\put(20,30){\line(1,1){10}} \put(40,30){\line(-1,1){10}}
\put(30,0){\makebox(0,0){Fig.5}}
\end{picture}
\end{center}

\bigskip

This point of view will be considered in a future work \cite{rappr}, in which
it will be applied to many classes of lattice paths.

\bigskip

Now we come to describe our construction of $\mathcal{D}_{n+1}$ starting
from $\mathcal{D}_{n}$. The first thing to observe is that it is possible
to determine a very natural chain partition of any Dyck lattice,
suggested by a particular method of enumeration, called the \emph{ECO
method}, which in this context reveals some unexpected peculiarities in
exploring order-theoretic properties of a combinatorial structure.

\bigskip

The ECO method was introduced in a series of articles, due to
Pinzani et al., in which they fruitfully apply a particular tool,
namely \emph{succession rules} (deeply studied by West and others,
see e.g. \cite{W}), to give a purely combinatorial construction
performing a local recursive expansion of a given class of
objects. Such a construction often allows to determine a
functional equation satisfied by the generating function of the
combinatorial structure under consideration. This, in turn, can be
solved (in many cases) to give an explicit expression for the
desired generating function. For a detailed survey concerning the
ECO method, we suggest \cite{eco}.

\bigskip

In this work, we are mainly interested in the ECO construction of Dyck
paths and in its properties with respect to the order introduced on Dyck
paths. So let us briefly recall how such construction works.

Consider a Dyck path $P$ of length $2n$, and suppose that the
length of its last descent (i.e., of its last sequence of fall
steps) is $k$. Then we construct $k+1$ Dyck paths of length $2n+2$
starting from $P$ (they will be called the \emph{sons} of $P$)
simply by inserting a peak (i.e., a rise step followed by a fall
step) in every point of its last descent. In fig. 6 it is shown in
a concrete example how this construction works.

\bigskip

\begin{center}
\setlength{\unitlength}{1.8mm}
\begin{picture}(120,25)
\put(0,17){\circle*{0.3}} \put(2,19){\circle*{0.3}}
\put(4,21){\circle*{0.3}} \put(6,19){\circle*{0.3}}
\put(8,21){\circle*{0.3}} \put(10,23){\circle*{0.3}}
\put(12,21){\circle*{0.3}} \put(14,19){\circle*{0.3}}
\put(16,17){\circle*{0.3}}
\put(20,19){\makebox(0,0){$\rightarrow$}}
\put(0,17){\line(1,1){2}} \put(2,19){\line(1,1){2}}
\put(4,21){\line(1,-1){2}} \put(6,19){\line(1,1){2}}
\put(8,21){\line(1,1){2}} \put(10,23){\line(1,-1){2}}
\put(12,21){\line(1,-1){2}} \put(14,19){\line(1,-1){2}}
\put(24,17){\circle*{0.3}} \put(26,19){\circle*{0.3}}
\put(28,21){\circle*{0.3}} \put(30,19){\circle*{0.3}}
\put(32,21){\circle*{0.3}} \put(34,23){\circle*{0.3}}
\put(36,21){\circle*{0.3}} \put(38,19){\circle*{0.3}}
\put(40,17){\circle*{0.3}} \put(42,19){\circle*{0.3}}
\put(44,17){\circle*{0.3}} \put(24,17){\line(1,1){2}}
\put(26,19){\line(1,1){2}} \put(28,21){\line(1,-1){2}}
\put(30,19){\line(1,1){2}} \put(32,21){\line(1,1){2}}
\put(34,23){\line(1,-1){2}} \put(36,21){\line(1,-1){2}}
\put(38,19){\line(1,-1){2}} \put(40,17){\line(1,1){2}}
\put(42,19){\line(1,-1){2}} \put(45,16){\makebox(0,0){,}}
\put(46,17){\circle*{0.3}} \put(48,19){\circle*{0.3}}
\put(50,21){\circle*{0.3}} \put(52,19){\circle*{0.3}}
\put(54,21){\circle*{0.3}} \put(56,23){\circle*{0.3}}
\put(58,21){\circle*{0.3}} \put(60,19){\circle*{0.3}}
\put(62,21){\circle*{0.3}} \put(64,19){\circle*{0.3}}
\put(66,17){\circle*{0.3}} \put(46,17){\line(1,1){2}}
\put(48,19){\line(1,1){2}} \put(50,21){\line(1,-1){2}}
\put(52,19){\line(1,1){2}} \put(54,21){\line(1,1){2}}
\put(56,23){\line(1,-1){2}} \put(58,21){\line(1,-1){2}}
\put(60,19){\line(1,1){2}} \put(62,21){\line(1,-1){2}}
\put(64,19){\line(1,-1){2}} \put(67,16){\makebox(0,0){,}}
\put(24,3){\circle*{0.3}} \put(26,5){\circle*{0.3}}
\put(28,7){\circle*{0.3}} \put(30,5){\circle*{0.3}}
\put(32,7){\circle*{0.3}} \put(34,9){\circle*{0.3}}
\put(36,7){\circle*{0.3}} \put(38,9){\circle*{0.3}}
\put(40,7){\circle*{0.3}} \put(42,5){\circle*{0.3}}
\put(44,3){\circle*{0.3}} \put(24,3){\line(1,1){2}}
\put(26,5){\line(1,1){2}} \put(28,7){\line(1,-1){2}}
\put(30,5){\line(1,1){2}} \put(32,7){\line(1,1){2}}
\put(34,9){\line(1,-1){2}} \put(36,7){\line(1,1){2}}
\put(38,9){\line(1,-1){2}} \put(40,7){\line(1,-1){2}}
\put(42,5){\line(1,-1){2}} \put(45,2){\makebox(0,0){,}}
\put(46,3){\circle*{0.3}} \put(48,5){\circle*{0.3}}
\put(50,7){\circle*{0.3}} \put(52,5){\circle*{0.3}}
\put(54,7){\circle*{0.3}} \put(56,9){\circle*{0.3}}
\put(58,11){\circle*{0.3}} \put(60,9){\circle*{0.3}}
\put(62,7){\circle*{0.3}} \put(64,5){\circle*{0.3}}
\put(66,3){\circle*{0.3}} \put(46,3){\line(1,1){2}}
\put(48,5){\line(1,1){2}} \put(50,7){\line(1,-1){2}}
\put(52,5){\line(1,1){2}} \put(54,7){\line(1,1){2}}
\put(56,9){\line(1,1){2}} \put(58,11){\line(1,-1){2}}
\put(60,9){\line(1,-1){2}} \put(62,7){\line(1,-1){2}}
\put(64,5){\line(1,-1){2}} \put(30,0){\makebox(0,0){Fig.6}}
\end{picture}
\end{center}

\bigskip

By performing this construction on all Dyck paths of length $2n$,
one obtains every Dyck path of length $2n+2$ exactly once.
Moreover, such a construction possesses strong recursive
properties, which can be encoded by the following succession rule:
\begin{equation}\label{regcat}
\Omega :\left\{
\begin{array}{l}
(2) \\
(k) \rightsquigarrow (2)(3)\ldots (k)(k+1).
\end{array}
\right.
\end{equation}

In the above expression it has to be intended that every Dyck path
whose last descent has length $k-1$ (encoded by the label $(k)$)
produces $k$ sons (through the above ECO construction) having last
descents of length $1,2,\ldots ,k-1,k$ respectively. A possible
way of visualizing this succession rule is to draw its
\emph{generating tree}, that is the infinite rooted labelled tree
in which the root is labelled $(2)$ and every node labelled $(k)$
has $(k)$ sons labelled $(2),(3),\ldots ,(k),(k+1)$ respectively.
Thus the number of Dyck paths of length $2n$ is precisely the
number of nodes at level $n$ of the generating tree. In this
situation, we say that the sequence of Catalan numbers is the
sequence \emph{defined by} (or related to, or associated with) the
rule (\ref{regcat}).

\bigskip

The above ECO construction of Dyck paths provides a partition of the sets of
Dyck paths of any fixed length. So one can consider the equivalence classes
of those Dyck paths having the same father. The following
result is stated without proof, since it is completely trivial; nevertheless,
it contains the main order-theoretic feature of the ECO construction of Dyck
paths.
\begin{teor}
The set of the sons of a Dyck path of length $2n-2$ forms a
saturated chain in the lattice $\mathcal{D}_{n}$ (i.e. a chain
such that for any two elements $x$ and $y$ in it, $x<y$ implies
that $x\preceq y$). Therefore, the ECO method provides a partition
into saturated chains of every Dyck lattice, which we will refer
to as the \emph{ECO-partition} of $\mathcal{D}_{n}$.
\end{teor}

In fig. 7 the ECO-partition of $\mathcal{D}_4$ is drawn.

\bigskip

\begin{center}
\setlength{\unitlength}{0.5mm}
\begin{picture}(40,70)
\put(50,10){\circle*{2}} \put(40,20){\circle*{2}}
\put(60,20){\circle*{2}} \put(50,30){\circle*{2}}
\put(60,40){\circle*{2}} \put(20,20){\circle*{2}}
\put(10,30){\circle*{2}} \put(30,30){\circle*{2}}
\put(20,40){\circle*{2}} \put(30,50){\circle*{2}} \put(
0,40){\circle*{2}} \put(10,50){\circle*{2}}
\put(20,60){\circle*{2}} \put(30,70){\circle*{2}}
\put(50,10){\line(1,1){10}} \put(40,20){\line(1,1){10}}
\put(50,30){\line(1,1){10}} \put(20,20){\line(1,1){10}}
\put(10,30){\line(1,1){10}} \put(20,40){\line(1,1){10}} \put(
0,40){\line(1,1){10}} \put(10,50){\line(1,1){10}}
\put(20,60){\line(1,1){10}} \multiput(50,10)(-2,2){5}{\circle*{1}}
\multiput(60,20)(-2,2){5}{\circle*{1}}
\multiput(20,20)(-2,2){5}{\circle*{1}}
\multiput(30,30)(-2,2){5}{\circle*{1}}
\multiput(10,30)(-2,2){5}{\circle*{1}}
\multiput(20,40)(-2,2){5}{\circle*{1}}
\multiput(30,50)(-2,2){5}{\circle*{1}}
\multiput(50,10)(-2,0.6){15}{\circle*{1}}
\multiput(40,20)(-2,0.6){15}{\circle*{1}}
\multiput(60,20)(-2,0.6){15}{\circle*{1}}
\multiput(50,30)(-2,0.6){15}{\circle*{1}}
\multiput(60,40)(-2,0.6){15}{\circle*{1}}
\put(20,0){\makebox(0,0){Fig.7}}
\end{picture}
\end{center}

\bigskip

Observe that the rule in (\ref{regcat}) is not the only possible rule defining
Catalan numbers; in \cite{equivalenze} many different rules for Catalan
numbers are derived, and the same thing can be done for any other numerical
sequence. However, it turns out that this is the "best" one
from an order-theoretical point of view, meaning that the nice order
structure of the equivalence classes gives a considerable help in the
description of Dyck lattices.

To understand the structure of Dyck lattices it remains now to investigate how
the saturated chains of the ECO-partition of $\mathcal{D}_{n}$ are linked
together. This is precisely what we are going to do in the next pages.

\bigskip

First of all we describe another partition of $\mathcal{D}_{n}$ into
sublattices (not necessarily chains) which follows in a very natural way from
some elementary geometric properties of Dyck paths. More precisely, the
sublattices of such a partition are certain filters of $\mathcal{D}_{n-1}$.

Denote by $\mathcal{D}_{nk}$ the set of all Dyck paths starting with
precisely $k$ rise steps (so that the $(k+1)$-st step is a fall step). In
symbols:
\begin{equation}
\mathcal{D}_{nk}=\{ f\in \mathcal{D}_{n}\; |\; \forall i\leq k, f(i)=i, f(k+1)=k-1\}.
\end{equation}

This definition is meaningful only when $k\leq n$. In particular,
$\mathcal{D}_{n0}=\emptyset$ and $\mathcal{D}_{nn}=\{
\mathbf{1}\}$. In fig. 8 it is shown how a path in
$\mathcal{D}_{n5}$ starts.

\bigskip

\begin{center}
\setlength{\unitlength}{0.5mm}
\begin{picture}(80,60)
\put(0,10){\circle*{1.5}} \put(10,20){\circle*{1.5}}
\put(20,30){\circle*{1.5}} \put(30,40){\circle*{1.5}}
\put(40,50){\circle*{1.5}} \put(50,60){\circle*{1.5}}
\put(60,50){\circle*{1.5}} \put(0,10){\line(1,1){10}}
\put(10,20){\line(1,1){10}} \put(20,30){\line(1,1){10}}
\put(30,40){\line(1,1){10}} \put(40,50){\line(1,1){10}}
\put(50,60){\line(1,-1){10}}
\multiput(60,50)(2,0){20}{\circle*{0.5}}
\put(40,0){\makebox(0,0){Fig.8}}
\end{picture}
\end{center}

\bigskip

\begin{prop}\label{dnkpartizione} For every $k\leq n$, $\mathcal{D}_{nk}$ is a sublattice of
$\mathcal{D}_{n}$. Moreover, the family of sublattices $\{
\mathcal{D}_{nk}\; |\; 0<k\leq n\}$ constitutes a partition of $\mathcal{D}_{n}$.
\end{prop}

\emph{Proof.} Observe that, if two paths starts with the same
sequence of steps, then this happens for their join and meet as
well: hence the first part of the proposition follows. The second
part is obvious, since the starting sequence of rise steps is
uniquely determined for every Dyck path. \cvd

\bigskip

In particular, for $\mathcal{D}_{n1}$ we have the following, nice result.

\begin{prop} $\mathcal{D}_{n1}\simeq \mathcal{D}_{n-1}$.
\end{prop}

\emph{Proof.} Consider the function mapping each element $f\in
\mathcal{D}_{n1}$ into the element of $\mathcal{D}_{n-1}$ obtained
by removing the first two steps of $f$: this is a bijection which
preserves joins and meets. \cvd

\bigskip

The above, very simple proposition states that, for every $n$,
$\mathcal{D}_{n}$ contains a copy of $\mathcal{D}_{n-1}$, which is
precisely the ideal of the elements of $\mathcal{D}_{n}$ starting
with a peak.

\bigskip

Now consider the sublattice $\mathcal{D}_{nk}$ of $\mathcal{D}_{n}$, for
$k>1$. The next proposition asserts that this lattice ``lives" also inside
$\mathcal{D}_{n-1}$ as a particular filter.

\begin{prop}\label{incoll1} Denote by $F_{(n-1)(k-1)}$ the subset of
$\mathcal{D}_{n-1}$ of all Dyck paths starting with $k-1$ rise
steps (without any further hypothesis on the $k$-th step). Then
$F_{(n-1)(k-1)}$ is a filter, and is isomorphic to
$\mathcal{D}_{nk}$.
\end{prop}

\emph{Proof.} The fact that $F_{(n-1)(k-1)}$ is a filter of
$\mathcal{D}_{n-1}$ is immediate (the argument is the same as that
of proposition \ref{dnkpartizione}). Now consider any path in
$\mathcal{D}_{nk}$ and delete its first peak (which, by
definition, occurs precisely after the first $k-1$ rise steps). In
this way we obtain a path of $\mathcal{D}_{n-1}$ starting with
$k-1$ rise steps, that is an element of $F_{(n-1)(k-1)}$. This
correspondence is a meet- and join-preserving bijection between
$\mathcal{D}_{nk}$ and $F_{(n-1)(k-1)}$, so the proof is complete.
\cvd

\bigskip

\underline{\emph{Remark}}. Notice that, as a filter of
$\mathcal{D}_{n-1}$, $F_{(n-1)(k-1)}$ is the principal filter
generated by the least Dyck path starting with a hill of height
$k-1$, in the sense that it contains precisely all the paths
greater or equal than the above mentioned one.

\bigskip

Next we show how the sublattices $\mathcal{D}_{nk}$ are linked together to
form the whole $\mathcal{D}_{n}$.

\begin{prop}\label{incoll2} Consider the function $\varphi
_{nk}:\mathcal{D}_{n(k+1)}\longrightarrow \mathcal{D}_{nk}$
defined as follows. Let $f\in \mathcal{D}_{n(k+1)}$, so that
$f(k+1)=k+1$, $f(k+2)=k$; we set by definition $\varphi
_{nk}(f)=g\in \mathcal{D}_{nk}$, where $f(i)=g(i), \forall i\neq
k+1$, and $g(k+1)=k-1$. Then $\varphi _{nk}$ is a lattice
monomorphism (i.e., injective homomorphism). Moreover, if
$g=\varphi _{nk}(f)$, then $f\prec g$ in $\mathcal{D}_{n}$.
\end{prop}

\emph{Proof.} First observe that $\varphi _{nk}$ is well-defined,
since $\varphi _{nk}(f)$ is an element of $\mathcal{D}_{nk}$. If
$\varphi_ {nk}(f_{1})=\varphi _{nk}(f_{2})$, then
$f_{1}(i)=f_{2}(i)$ at least for every $i\neq k+1$; moreover
$f_{1}(k+1)=k+1=f_{2}(k+1)$ (since $f_{1},f_{2}\in
\mathcal{D}_{n(k+1)}$), so that $f_{1}\equiv f_{2}$. Now consider
$f,f'\in \mathcal{D}_{n(k+1)}$, so that $f$ and $f'$ coincide for
$i\leq k+2$. Then also $f\lor f'$ coincides with both $f$ and $f'$
for $i\leq k+2$. Therefore, for $i\neq k+1$, we have $\varphi
_{nk}(f\lor f')(i)=(f\lor f')(i)=(\varphi _{nk}(f)\lor \varphi
_{nk}(f'))(i)$, whereas $\varphi _{nk}(f\lor f')(k+1)=k-1=(\varphi
_{nk}(f)\lor \varphi _{nk}(f'))(k+1)$, and so $\varphi _{nk}(f\lor
f')=\varphi _{nk}(f)\lor \varphi _{nk}(f')$. The same argument can
be used also for meet. Finally, we have that $f\prec \varphi
_{nk}(f)$, since $\alpha (f)+2=\alpha (\varphi _{nk}(f))$ and
$f<g$ (recall that $\alpha (f)$ denotes the area of $f$, as
defined in proposition \ref{rango}). \cvd

\bigskip

Fig. 9 illustrates how the map $\varphi _{n3}$ works.

\bigskip

\begin{center}
\setlength{\unitlength}{0.5mm}
\begin{picture}(180,60)
\put(0,10){\circle*{1.5}} \put(10,20){\circle*{1.5}}
\put(20,30){\circle*{1.5}} \put(30,40){\circle*{1.5}}
\put(40,50){\circle*{1.5}} \put(50,40){\circle*{1.5}}
\put(0,10){\line(1,1){10}} \put(10,20){\line(1,1){10}}
\put(20,30){\line(1,1){10}} \put(30,40){\line(1,1){10}}
\put(40,50){\line(1,-1){10}}
\multiput(50,40)(2,0){10}{\circle*{0.5}}
\put(80,40){\makebox(0,0){$\longrightarrow$}}
\put(100,10){\circle*{1.5}} \put(110,20){\circle*{1.5}}
\put(120,30){\circle*{1.5}} \put(130,40){\circle*{1.5}}
\put(140,30){\circle*{1.5}} \put(150,40){\circle*{1.5}}
\put(100,10){\line(1,1){10}} \put(110,20){\line(1,1){10}}
\put(120,30){\line(1,1){10}} \put(130,40){\line(1,-1){10}}
\put(140,30){\line(1,1){10}}
\multiput(150,40)(2,0){10}{\circle*{0.5}}
\put(90,0){\makebox(0,0){Fig.9}}
\end{picture}
\end{center}

\bigskip

\underline{\emph{Remark}}. Observe that $\varphi
_{nk}(\mathcal{D}_{n(k+1)})$ is a filter of $\mathcal{D}_{nk}$.
More precisely, it is the filter of all the paths of
$\mathcal{D}_{nk}$ whose $(k+2)$-nd step is a rise step (the
reader can check that this is actually a filter of
$\mathcal{D}_{nk}$).

\bigskip

Let us summarize what we have done until now. We have discovered a
special partition of $\mathcal{D}_{n}$ into sublattices (the
$\mathcal{D}_{nk}$'s). We have shown that a certain filter of
$\mathcal{D}_{n-1}$ is isomorphic to $\mathcal{D}_{nk}$ \footnote{In
particular, $\mathcal{D}_{n1}$ is isomorphic to the whole
$\mathcal{D}_{n-1}$.}, which, in turn, is isomorphic to a filter of
$\mathcal{D}_{n(k-1)}$. As a consequence of this fact, we have that
$\mathcal{D}_{n}$ can be constructed starting from $\mathcal{D}_{n-1}$ or,
which is the same, from $\mathcal{D}_{n1}$ by glueing together some of its
filters. Notice that this is all we need to draw the Hasse diagram of
$\mathcal{D}_{n}$, since $x\prec y$ in $\mathcal{D}_{n}$ if and only if
either $x\prec y$ in $\mathcal{D}_{nk}$, for some $k$, or
$x\in \mathcal{D}_{nk}$, $y\in \mathcal{D}_{n(k+1)}$ and $y=\varphi
_{nk}(x)$.

\bigskip

Now, to give a complete formalization of our construction we need to define
a particular operation on lattices.

Given a lattice $\mathbf{L}$, let $F$ be a filter of $\mathbf{L}$. We define
the \emph{F-filtered doubling} of $\mathbf{L}$ to be the lattice
\begin{equation}
\mathbf{L}\times_{F}\mathbf{2}=\{ (x,n)\in L\times \{ 0,1\} \; |\;
n=1\Rightarrow x\in F\}
\end{equation}
endowed with coordinatewise meet and join.

\bigskip

\underline{\emph{Remark}}. We wish to point out that this
construction is not new at all. In fact, this is merely a
particular instance of a \emph{subdirect product} of lattices.
Recall that a subdirect product of two lattices $\mathbf{L}$ and
$\mathbf{M}$ is, by definition, any sublattice $\mathbf{N}$ of the
direct product $\mathbf{L}\times \mathbf{M}$ such that both the
projections of $\mathbf{N}$ onto $\mathbf{L}$ and $\mathbf{M}$ are
surjective homomorphisms (see, for example, \cite{CD}). So, an
F-filtered doubling of $\mathbf{L}$ is nothing more than a
particular subdirect product of the lattices $\mathbf{L}$ and
$\mathbf{2}=[\{ 0,1\} ;\lor ,\land ]$. By the way, this also
proves that $\mathbf{L}\times_{F}\mathbf{2}$ is actually a
lattice. Observe that, if $F=L$, then $\mathbf{L}\times
_{F}\mathbf{2}=\mathbf{L}\times \mathbf{2}$.

\bigskip

Now consider, for every positive integer $k<n$, all the filters $F_{(n-1)k}$
of $\mathcal{D}_{n-1}$ as defined above. So, for example, $F_{(n-1)1}$ is the
whole lattice $\mathcal{D}_{n-1}$ and $F_{(n-1)(n-1)}$ is the singleton of
the maximum of $\mathcal{D}_{n-1}$. We introduce the following recursive
notation:
\begin{eqnarray*}
\beta _{1}(\mathcal{D}_{n-1})=\mathcal{D}_{n-1}\times
_{F_{(n-1)1}}\mathbf{2}=\mathcal{D}_{n-1}\times \mathbf{2},
\\ \beta _{k}(\mathcal{D}_{n-1})=\beta _{k-1}(\mathcal{D}_{n-1})\times
_{F_{(n-1)k}}\mathbf{2},\qquad \textnormal{for $k<n$}.
\end{eqnarray*}

Observe that the operations performed make sense since, in general,
$F_{(n-1)k}$ is a filter of $\beta _{k-1}(\mathcal{D}_{n-1})$ (this is a
consequence of propositions \ref{incoll1} and \ref{incoll2}).

The following theorem condenses in a single formula the whole construction
of $\mathcal{D}_{n}$ from $\mathcal{D}_{n-1}$ described in this section,
making use of the notations introduced above.

\begin{teor} $\mathcal{D}_{n}=\beta _{n-1}(\mathcal{D}_{n-1})$.
\end{teor}

\subsection{Schr\"oder lattices}

As we have said in section 2, Schr\"oder paths do not fall within
the class of paths captured by our definition. Recall that a
Schr\"oder path is a lattice path starting from the origin, ending
on the $x$-axis, never falling below the $x$-axis and using steps
of the type $(1,1)$ (rise), $(1,-1)$ (fall) and $(2,0)$ (double
horizontal). Nevertheless, they can be naturally ordered in the
same way, and it can be shown that they constitute a lattice. It
turns out that the study of Schr\"oder lattices follows
essentially the same lines as those for Dyck lattices. Also in
this case, there is a privileged rule among all those defining
Schr\"oder numbers, which is precisely
\begin{equation}\label{regscr}
\Omega :\left\{
\begin{array}{l}
(2) \\
(2k) \rightsquigarrow (2)(4)^{2}(6)^{2}\ldots (2k-2)^{2}(2k)^{2}(2k+2).
\end{array}
\right.
\end{equation}

The ECO-construction of Schr\"oder paths encoded by the above rule
provides a partition of any Schr\"oder lattice $\mathcal{S}_{n}$
into saturated chain. Besides, the construction of
$\mathcal{S}_{n}$ starting from $\mathcal{S}_{n-1}$ can be carried
out in an analogous way to Dyck lattices: the interested reader
will find no difficulties in reproducing it. In fig. 10 we give
the Hasse diagram of $\mathcal{S}_{3}$ (lattice of Schr\"oder
paths of length 6).

\bigskip

\begin{center}
\setlength{\unitlength}{1mm}
\begin{picture}(80,100)
\put(60,10){\circle*{1}} \put(70,20){\circle*{1}}
\put(50,20){\circle*{1}} \put(60,30){\circle*{1}}
\put(70,40){\circle*{1}} \put(80,50){\circle*{1}}
\put(30,20){\circle*{1}} \put(40,30){\circle*{1}}
\put(20,30){\circle*{1}} \put(30,40){\circle*{1}}
\put(40,50){\circle*{1}} \put(50,60){\circle*{1}}
\put(10,40){\circle*{1}} \put(20,50){\circle*{1}}
\put(30,60){\circle*{1}} \put(40,70){\circle*{1}}
\put(0,50){\circle*{1}} \put(10,60){\circle*{1}}
\put(20,70){\circle*{1}} \put(30,80){\circle*{1}}
\put(40,90){\circle*{1}} \put(50,100){\circle*{1}}
\put(60,10){\line(1,1){10}} \put(50,20){\line(1,1){10}}
\put(60,30){\line(1,1){10}} \put(70,40){\line(1,1){10}}
\put(30,20){\line(1,1){10}} \put(20,30){\line(1,1){10}}
\put(30,40){\line(1,1){10}} \put(40,50){\line(1,1){10}}
\put(10,40){\line(1,1){10}} \put(20,50){\line(1,1){10}}
\put(30,60){\line(1,1){10}} \put(0,50){\line(1,1){10}}
\put(10,60){\line(1,1){10}} \put(20,70){\line(1,1){10}}
\put(30,80){\line(1,1){10}} \put(40,90){\line(1,1){10}}
\put(60,10){\line(-1,1){10}} \put(70,20){\line(-1,1){10}}
\put(30,20){\line(-1,1){10}} \put(40,30){\line(-1,1){10}}
\put(20,30){\line(-1,1){10}} \put(30,40){\line(-1,1){10}}
\put(40,50){\line(-1,1){10}} \put(50,60){\line(-1,1){10}}
\put(10,40){\line(-1,1){10}} \put(20,50){\line(-1,1){10}}
\put(30,60){\line(-1,1){10}} \put(40,70){\line(-1,1){10}}
\put(60,10){\line(-3,1){30}} \put(70,20){\line(-3,1){30}}
\put(50,20){\line(-3,1){30}} \put(60,30){\line(-3,1){30}}
\put(70,40){\line(-3,1){30}} \put(80,50){\line(-3,1){30}}
\put(20,10){\makebox(0,0){$\mathcal{S}_{3}$}}
\put(40,0){\makebox(0,0){Fig.10}}
\end{picture}
\end{center}

\section{Conclusions and open problems}

In this paper we would like to begin a deep investigation on the
order-theoretic aspects of the ECO method, so it should be
intended only as a first step in this direction. We hope that the
present work will be followed soon by analogous ones concerning
several other classes of combinatorial objects naturally definable
by means of some ECO-construction. As far as paths are concerned,
the following problems seems to be particularly challenging.

\begin{enumerate}
\item Recall that a finite sequence of numbers $a_{0},\ldots a_{n}$
is said to be \emph{unimodal} whenever there exists an index $k$ such
that $a_{0}\leq a_{1}\leq \ldots \leq a_{k}$ and $a_{k}\geq a_{k+1}
\geq \ldots \geq a_{n}$. It seems natural to think that Dyck lattices are
rank-unimodal (i.e. the sequence given by the elements of the same rank
is unimodal), even if we have not succeeded in proving it yet. This problem
was also mentioned in \cite{unimod}, both for Dyck and
Schr\"oder lattices, where it is presented as an ``intriguing open
question". Obviously the same question can be asked for any other lattice
of paths. In
\cite{momenti} the author studies problems related with the area and the
moments of Dyck paths, and in \cite{libro}, exercise 5.2.12 gives formulas
for the area under various paths: these references could be of some help in
tackling this problem.
\item The approach used in the study of Dyck and Schr\"oder lattices
does not work for Motzkin lattices. In particular, the ECO decomposition for
Motzkin lattices does not provide a good description, as it happens for Dyck
and Schr\"oder lattices. In fact, it is still true that the sons of
a Motzkin path constitute a chain of the Motzkin lattice they belong to, but
such a chain is not saturated (see \cite{eco} for a precise description of
the usual ECO construction of Motzkin paths). Therefore we are not helped
by the ECO method
in drawing the Hasse diagrams of Motzkin lattices. So the study of such
lattices (as well as of other lattices of paths of some interest in
combinatorics) should be done, maybe using different techniques from those
developed in the present work.
\end{enumerate}

{\footnotesize \textbf{Acknowledgments.} The authors wish to thank
Sri Gopal Mohanty for drawing their attention to the references
\cite{nara,nara2}. They also thank the referees for useful
suggestions and comments.}

\end{document}